\documentclass[11pt,reqno]{amsart}
\usepackage{amssymb,amsmath,epsfig}
\usepackage{amsfonts}
\newtheorem{theorem}{Theorem}

\newtheorem{lemma}[theorem]{Lemma}
\newtheorem{corollary}[theorem]{Corollary}

\newtheorem*{theorem*}{Theorem}
\title{AN ITERATIVE-BIJECTIVE APPROACH TO GENERALIZATIONS OF SCHUR's THEOREM}

\date{\today}
\author{Sylvie Corteel and Jeremy Lovejoy}
\address{CNRS, PRISM, Universit\'e de Versailles Saint Quentin,
45 Avenue des Etats Unis, 78035 Versailles Cedex, FRANCE}
\address{CNRS, LABRI, Universit\'e Bordeaux I,
351 Cours de la lib\'eration, 33405 Talence Cedex, FRANCE}
\email{syl@prism.uvsq.fr} \email{lovejoy@labri.fr}
\thanks{The second author was partially supported by the European Commission's IHRP
  Programme, grant HPRN-CT-2001-00272, ``Algebraic Combinatorics in
  Europe''}
\begin{document}
\begin{abstract}
We start with a bijective proof of Schur's theorem due to Alladi
and Gordon and describe how a particular iteration of it leads to
some very general theorems on colored partitions.  These theorems
imply a number of important results, including Schur's theorem,
Bressoud's generalization of a theorem of G\"ollnitz, two of
Andrews' generalizations of Schur's theorem, and the
Andrews-Olsson identities.
\end{abstract}
\maketitle
\section{Introduction}
In $1926$ Schur \cite{Sc1} proved that the number of partitions of
$m$ into distinct parts not divisible by $3$ is equal to the
number of partitions of $m$ where parts differ by at least $3$ and
multiples of $3$ differ by at least $6$.  Over the years there
have been a number of proofs of this theorem (e.g.
\cite{Al-Go1,An.5,An1.5,An4,Be1,Br1}), including a few
delightfully simple combinatorial arguments \cite{Al-Go1,Be1,Br1}.
In \cite{Al-Go1}, Alladi and Gordon showed how to deduce Schur's
theorem from an interpretation of the infinite product
\begin{equation} \label{twoproduct}
\prod_{k=1}^{\infty} (1+y_1q^k)(1+y_2q^k)
\end{equation}
as a generating function for partitions whose parts come in $3$
colors:

\begin{theorem}[Alladi-Gordon] \label{alladi}
The number of pairs of partitions $(\mu_1,\mu_2)$, where $\mu_r$
is a partition into $x_r$ distinct parts and the sum of all of the
parts is $m$, is equal to the number of partitions of $m$ into
distinct parts occurring in $3$ colors (labelled $1,2$, and $3$)
such that (i) the part $1$ does not occur in color $3$, (ii)
consecutive parts differ by at least $2$ if the larger has color
$3$ or if the larger has color $1$ and the smaller has color $2$,
(iii) $x_j$ is the number of parts with color $j$ plus the number
of parts with color $3$.
\end{theorem}

Among their proofs of this theorem is an attractive bijective
argument which adapts some ideas of Bressoud \cite{Br.5}. In this
paper we describe a particular iteration of this bijection, which
leads an interpretation of the infinite product
\begin{equation} \label{product}
\prod_{k=1}^{\infty} (1+y_1q^k)(1+y_2q^k)\cdots (1+y_nq^k)
\end{equation}
as a generating function for certain partitions where the parts
come in $2^n-1$ colors.  To state the theorems, we require some
notation. For $t$-colored partitions, we denote the $t$ colors by
the natural numbers $(1,...,t)$, with the parts ordered first
according to size and then according to color. We write
$\omega(c)$ for the number of powers of $2$ occurring in the
binary representation of $c$, and $v(c)$ (resp. $z(c)$) for the
smallest (resp. largest) power of $2$ occurring in this
representation.  The function $\delta(c,d)$ is equal to $1$ if
$z(c) < v(d)$ and is $0$ otherwise.  Finally, we use $c_i$ to
denote the color of a part $\lambda_i$.

\begin{theorem} \label{lovejoy1}
Let $A(x_1,\ldots,x_n;m)$ denote the number of $n$-tuples
$(\mu_1,\mu_2,...,\mu_n)$, where $\mu_r$ is a partition into $x_r$
distinct parts, and the sum of all of the parts is $m$. Let
$B(x_1,\ldots,x_n;m)$ denote the number of partitions $\lambda_1 +
\cdots +\lambda_s$ of $m$ into distinct parts occurring in $2^n -
1$ colors, where (i) $\lambda_s \geq \omega(c_s)$, (ii) $x_r$ of
the colors $c_i$ have $2^{r-1}$ in their binary representations,
and (iii) $\lambda_i - \lambda_{i+1} \geq \omega(c_i) +
\delta(c_i,c_{i+1})$. Then $A(x_1,\ldots,x_n;m) =
B(x_1,\ldots,x_n;m).$
\end{theorem}

For example, using the notation $(\lambda_1,...,\lambda_s)$ to
represent the sum $\lambda_1 + \cdots + \lambda_s$ and the symbol
$\epsilon$ to denote the empty partition, the $4$-tuples of
partitions counted by $A(2,0,1,1;9)$ are
$$
\begin{gathered}
((6,1),\epsilon,1,1),((5,2),\epsilon,1,1),((4,3),\epsilon,1,1),((5,1),\epsilon,2,1),((5,1),\epsilon,1,2),
((4,2),\epsilon,2,1),
\\
((4,2),\epsilon,1,2),((4,1),\epsilon,3,1),((4,1),\epsilon,2,2),((4,1),\epsilon,1,3)
((3,2),\epsilon,3,1),((3,2),\epsilon,2,2),
\\
((3,2),\epsilon,1,3)((3,1),\epsilon,4,1),((3,1),\epsilon,3,2),
((3,1),\epsilon,2,3),((3,1),\epsilon,1,4),((2,1),\epsilon,5,1),
\\
((2,1),\epsilon,4,2),((2,1),\epsilon,3,3),((2,1),\epsilon,2,4),((2,1),\epsilon,1,5)
\end{gathered}
$$
and the partitions counted by $B(2,0,1,1;9)$ are
$$
\begin{gathered}
(8_{13},1_1),(7_{13},2_1),(7_9,2_5),(7_5,2_9),(6_{13},3_1),(6_{9},3_5),(6_5,3_{9}),
(6_1,3_{13}),(5_1,4_{13}), \\
(6_9,2_4,1_1),(6_{12},2_1,1_1),(6_5,2_8,1_1),(5_9,3_4,1_1),(5_4,3_9,1_1),(5_8,3_5,1_1),
(5_9,3_1,1_4),(5_1,3_9,1_4),
\\
(5_5,3_1,1_8),(5_{12},3_1,1_1),(5_1,3_{12},1_1),(4_8,3_1,2_5),(4_4,3_1,2_9).
\end{gathered}
$$

Actually, a closer look at the proof of Theorem \ref{lovejoy1}
will reveal that it can be extended by letting $\mu_1$ be either a
partition into distinct parts congruent to $R$ modulo $M$ or a
partition into parts that differ by at least $M$. These cases
correspond to the products
\begin{equation} \prod_{k=1}^{\infty}(1+y_1q^{(k-1)M
+R})(1+y_2q^{k})\cdots (1+y_nq^{k}) \label{product1}
\end{equation}
and
\begin{equation}
\sum_{k=0}^{\infty} \frac{y_1^kq^{M(k(k-1)/2)
+k}}{(1-q)(1-q^2)\cdots(1-q^k)}\prod_{k=1}^{\infty}(1+y_2q^{k})\cdots
(1+y_nq^{k}). \label{product2}
\end{equation}
Let $\omega_e(c)$ denote the number of even powers of $2$ in the
binary representation of $c$.

\begin{theorem} \label{lovejoy2}
Let $A_{R,M}(x_1,\ldots,x_n;m)$ denote the number of $n$-tuples
$(\mu_1,\mu_2,...,\mu_n)$, where each $\mu_r$ is a partition into
$x_r$ distinct parts, $\mu_1$ is a partition into distinct parts
congruent to $R$ modulo $M$, and the sum of all of the parts is
$m$.  Let $B_{R,M}(x_1,\ldots,x_n;m)$ denote the number of
partitions $\lambda_1 + \cdots + \lambda_s$ of $m$ counted by
$B(x_1,\ldots,x_n;m)$ such that (i) if $\lambda_{S}$ denotes the
smallest part with odd color, then $\lambda_{S} \equiv R +
\sum_{\ell = S}^{s} \omega_e(c_{\ell}) \pmod{M}$ and (ii) if
$\lambda_i \geq \lambda_j$ are any two parts with odd color, then
$\lambda_i \equiv \lambda_j + \sum_{\ell = i}^{j-1}
\omega_e(\lambda_{\ell}) \pmod{M}$.  Then
$A_{R,M}(x_1,\ldots,x_n;m) = B_{R,M}(x_1,\ldots,x_n;m).$
\end{theorem}

\begin{theorem} \label{lovejoy3}
Let $A_{M}(x_1,\ldots,x_n;m)$ denote the number of $n$-tuples
$(\mu_1,\mu_2,...,\mu_n)$, where each $\mu_r$ is a partition into
$x_r$ distinct parts and $\mu_1$ is a partition into parts
differing by at least $M$, and the sum of all of the parts is $m$.
Let $B_{M}(x_1,\ldots,x_n;m)$ denote the number of partitions
$\lambda_1 + \cdots + \lambda_s$ of $m$ counted by
$B(x_1,\ldots,x_n;m)$ such that if $\lambda_i \geq \lambda_j$ are
any two parts with odd color, then $\lambda_i - \lambda_j \geq M +
\sum_{\ell = i}^{j-1}\omega_e(\lambda_{\ell})$. Then
$A_{M}(x_1,\ldots,x_n;m) = B_{M}(x_1,\ldots,x_n;m).$
\end{theorem}

Theorems \ref{lovejoy1} - \ref{lovejoy3} are closely related to a
number of important results in the theory of partitions.  For
example, by appropriately defining a conjugation on the partitions
counted by $B(x_1,\ldots,x_n;m)$ we will arrive at Theorem
\ref{lovejoy1conj} below, which is a generalization of the
Andrews-Olsson identities \cite{An-Ol1} and which was proven by
Bessenrodt \cite[Theorem 2.4, $\mathcal{C'} = \emptyset$]{Be2} and
stated by Alladi \cite[Theorem 15]{Al1}. Here we shall use
uncolored parts as well as colored parts, assuming that an
uncolored part of a given size occurs before all other parts of
that size.

\begin{theorem} \label{lovejoy1conj}
Let $C(x_1,\ldots,x_n;m)$ denote the number of partitions
$\lambda_1 + \cdots +\lambda_s$ of $m$ into distinct parts
occurring either in $n$ colors or uncolored, where (i) the
smallest part is colored, (ii) $x_r$ of the parts have color $r$,
and (iii) $\lambda_i - \lambda_{i+1} \leq 1$, with strict
inequality if $c_{i+1} < c_i$ or $\lambda_i$ is uncolored.  Then
$A(x_1,\ldots,x_n;m) = B(x_1,\ldots,x_n;m) = C(x_1,\ldots,x_n;m).$
\end{theorem}

Theorem \ref{lovejoy1conj} and some of the many other partition
theorems contained in Theorems \ref{lovejoy1}- \ref{lovejoy3} are
discussed in more detail in Section 5.  In the following section
we review the Alladi-Gordon bijective proof of Theorem
\ref{alladi} and explain our proof of Theorem \ref{lovejoy1} in
the case $n=3$. In Section 3 we undertake the proof of Theorem
\ref{lovejoy1} in full generality and in Section 4 we prove the
extensions, Theorems \ref{lovejoy2} and \ref{lovejoy3}.

\section{Two basic cases}
Although the basic idea behind the proofs of Theorems
\ref{lovejoy1} - \ref{lovejoy3} is a simple one, the amount of
notation required may obscure this fact. Therefore, we shall
present the cases $n=2$ and $n=3$ in detail. First, we review the
Alladi-Gordon bijective proof of Theorem \ref{alladi}, which is
the case $A(x_1,x_2;m) = B(x_1,x_2;m)$ of Theorem \ref{lovejoy1}.
We begin with a partition $\lambda$ into $x_1$ distinct parts
colored by $1$ and a partition $\tau$ into $x_2$ distinct parts
colored by $2$.

$\bullet$ {\bf Step 1.} For each part $k$ of $\tau$ that is less
than or equal to the number of parts of $\lambda$, we add $1$ to
the first $k$ parts of $\lambda$ and $2$ to the color of
$\lambda_k$. We then have the difference conditions
$$
\lambda_i-\lambda_{i+1}\ge \omega(c_i),
$$
for $c_i,c_{i+1} \neq 2$.  Notice that all parts with color $3$
are bigger than $1$.

$\bullet$ {\bf Step 2.} Now write the unused parts of $\tau$ in
decreasing order to the left of the parts from $\lambda$. Remove a
staircase, i.e, subtract 0 from the smallest part, 1 from the next
smallest, and so on. We therefore get the difference conditions
$$
\lambda_i-\lambda_{i+1}\ge \omega(c_i)-1
$$
for $c_i,c_{i+1} \neq 2$.

$\bullet$ {\bf Step 3.} Each part $\tau_j$ of color $2$ remaining
in $\tau$ is inserted in $\lambda$ after the smallest part that is
bigger than $\tau_j$.  We now have the difference conditions

\begin{equation*}
\lambda_i-\lambda_{i+1} \ge
\begin{cases}
\omega(c_i)-1, & c_i,c_{i+1}=1,3, \\
0, & c_i=2, \\
1, & c_i=1,3,\ c_{i+1}=2.
\end{cases}
\end{equation*}

$\bullet$ {\bf Step 4.} In each case above, the minimum difference
is exactly $\omega(c_i) + \delta({c_i,c_{i+1}})-1$. We add back
the staircase removed in Step 2 and we have
$$\lambda_i-\lambda_{i+1}\ge \omega(c_i)+\delta({c_i,c_{i+1}}).$$
This is condition $(ii)$ of Theorem \ref{alladi}.  Conditions
$(i)$ and $(iii)$ are straightforward.

\medskip

For example, starting with $\lambda=(8_1,3_1,2_1,1_1)$ and
$\tau=(10_2,5_2,3_2,2_2)$, we perform the steps of the bijection:
\begin{eqnarray*}
(\tau,\lambda) &\iff & ((10_2,5_2),(10_1,5_3,3_3,1_1)) \hskip.5 in
\text{(Step 1)}
\\ &\iff& ((5_2,1_2),(7_1,3_3,2_3,1_1)) \hskip.5in \text{(Step 2)} \\
&\iff& (\epsilon,(7_1,5_2,3_3,2_3,1_2,1_1)) \hskip.5in \text{(Step 3)}\\
&\iff& (12_1,9_2,6_3,4_3,2_2,1_1) \hskip.5in \text{(Step 4)}.
\end{eqnarray*}

Now, since the result of the above process is another partition
into distinct parts, it is natural to attempt to apply the
bijection again, starting with the partition $\lambda$ into
distinct parts having 3 colors satisfying $\lambda_i-\lambda_{i+1}
\ge \omega(c_i)+\delta({c_i,c_{i+1}})$ and $\lambda_i\ge
\omega({c_i})$ and a new partition, $\tau$, into distinct parts
occurring in the color $4$.  Let us see what happens when we try
to repeat the steps above.

$\bullet$ {\bf Step 1.} For each part $k$ of $\tau$ that is less
than or equal to the number of parts of $\lambda$, we add $1$ to
the first $k$ parts of $\lambda$ and $4$ to the color of
$\lambda_k$. We then have the difference conditions
$$
\lambda_i-\lambda_{i+1}\ge \omega(c_i) + \delta^*({c_i,c_{i+1}}),
$$
and $\lambda_i\ge \omega({c_i})$ for $c_i,c_{i+1}\neq 4$. Here
$\delta^*({c_i,c_{i+1}})=\delta({c_i,c_{i+1}})$ if $c_i<4$ and
$\delta({c_i-4,c_{i+1}})$ otherwise.

$\bullet$ {\bf Step 2.} Now write the unused parts of $\tau$ to
the left of the parts from $\lambda$. Remove a staircase, i.e
subtract 0 from the smallest part, 1 from the next smallest, and
so on. We get
$$
\lambda_i-\lambda_{i+1}\ge \omega(c_i)+\delta^*({c_i,c_{i+1}})-1,
$$
for $c_i,c_{i+1}\neq 4$.

$\bullet$ {\bf Step 3.} Now we take largest part $\tau_1$
remaining in $\tau$ and insert it into $\lambda$ after the
smallest part that is bigger than $\tau_1$.
It is easy to check that after such an insertion, we have
\begin{equation} \label{badcond}
\lambda_i-\lambda_{i+1} \ge
\begin{cases}
\omega(c_i)-\delta^*({c_i,c_{i+1}})-1, & c_i,c_{i+1} \neq 4 \\
0 = \omega(c_i)-\delta({c_i,c_{i+1}})-1, & c_i=4 \\
1 = \omega(c_i)-\delta({c_i,c_{i+1}})-u(c_i)-1, & c_i \neq 4,\
c_{i+1}=4.
\end{cases}
\end{equation}
Here $u(j)=1$, if $j=3$ or $7$, and $u(j) = 0$ otherwise.

Unlike the situation for $n=2$, we do not always have the
condition
\begin{equation} \label{goodcond}
\lambda_i-\lambda_{i+1}\ge \omega(c_i)+\delta({c_i,c_{i+1}})-1.
\end{equation}
There are four cases where \eqref{badcond} and \eqref{goodcond} do
not match up: $(i)$ If $c_i=7$ and $c_{i+1}=4$, the minimal
difference between $\lambda_i$ and $\lambda_{i+1}$ is
$1=\omega(c_i) + \delta({c_i,c_{i+1}})-2$, $(ii)$ When $c_i=5$ and
$c_{i+1}=6$, this difference is $2=\omega(c_i) +
\delta({c_i,c_{i+1}})$, $(iii)$  For $c_i=3$ and $c_{i+1}=4$, it
is $1=\omega(c_i) + \delta({c_i,c_{i+1}})-2$, and $(iv)$ for
$c_i=5$ and $c_{i+1}=2$, it is
$2=\omega(c_i)+\delta({c_i,c_{i+1}})$.

In the second and fourth cases, the guaranteed minimum difference
between $\lambda_i$ and $\lambda_{i+1}$ is too large. In the first
and third cases, it is too small.  For the latter cases, we shall
remedy the problem by redistributing the powers of $2$ occurring
in the colors $c_i$ and $c_{i+1}$. Specifically, if $c_i=7$,
$c_{i+1}=4$, and $\lambda_i-\lambda_{i+1}=1$, then let $c_i=5$ and
$c_{i+1}=6$. If $c_i=3$, $c_{i+1}=4$, and
$\lambda_i-\lambda_{i+1}=1$, then let $c_i=5$ and $c_{i+1}=2$.
Notice that the new colors correspond exactly to those in cases
$(ii)$ and $(iv)$ above, with the difference between $\lambda_i$
and $\lambda_{i+1}$ exactly one less than the minimum difference
in \eqref{badcond} corresponding to these two cases.  This is a
double bonus. First, it makes the change of colors described above
bijective, and second, it makes the minimum difference what we
want for the theorem.

We repeat Step 3 with the largest part remaining in $\tau$ and
continue until all the parts of color $4$ are inserted in
$\lambda$.  We then have \eqref{goodcond} for all $i$.

$\bullet$ {\bf Step 4.} Finally, we can add back the staircase and
we have
$$\lambda_i-\lambda_{i+1}\ge \omega(c_i)+\delta({c_i,c_{i+1}}).$$
This is condition $(iii)$ of Theorem \ref{lovejoy1}.  Conditions
$(i)$ and $(ii)$ are again straightforward.

\medskip

For example, we start with $\lambda= (12_1,9_2,6_3,4_3,2_2,1_1)$
and $\tau=(17_4,11_4,8_4,6_4,3_4,1_4)$, following the steps of the
bijection:
\begin{eqnarray*}
(\tau,\lambda) &\iff &
((17_4,11_4,8_4),(15_5,11_2,8_7,5_3,3_2,2_5)) \hskip.5 in
\text{(Step 1)}
\\ &\iff&  ((9_4, 4_4, 2_4),(10_5,7_2,5_7,3_3,2_2,2_5)) \hskip.5in \text{(Step 2)} \\
&\iff&  ((4_4, 2_4),(10_5,9_4,7_2,5_7,3_3,2_2,2_5)) \hskip.5in \text{(Step 3)} \\
&\iff&  ((2_4),(10_5,9_4,7_2,\underbrace{5_7,4_4},3_3,2_2,2_5))  \\
&\iff&  ((2_4),(10_5,9_4,7_2,{5_5,4_6},3_3,2_2,2_5))  \\
&\iff&
(\epsilon,(10_5,9_4,7_2,{5_5,4_6},\underbrace{3_3,2_4},2_2,2_5)) \\
&\iff& (\epsilon,(10_5,9_4,7_2,{5_5,4_6},{3_5,2_2},2_2,2_5)) \\
&\iff& (18_5,16_4,13_2,{10_5,8_6},{6_5,4_2},3_2,2_5) \hskip.5in
\text{(Step 4)}
\end{eqnarray*}
Here the underbraces indicate where a reassignment of colors needs
to take place.\\

The proof of Theorem \ref{lovejoy1} is an iteration of the above
process, described for a general $n$ in the following section.

\section{Proof of Theorem \ref{lovejoy1}}
\setcounter{equation}{0} \setcounter{theorem}{0}

The proof is by induction. For $n=1$ the theorem is a tautology.
Now suppose that it is true for a given natural number $n-1$ and
that we can map any partition counted by
$A(x_1,x_2,...,x_{n-1};m)$ to a partition counted by
$B(x_1,x_2,...,x_{n-1};m)$.  Take a partition counted by
$A(x_1,\ldots,x_n;m)$ and break it into two partitions: the first
is a partition counted by $A(x_1,x_2,...,x_{n-1};m')$ and the
second is a partition of $m-m'$ into $x_n$ distinct parts. We
apply the map to the first partition to get a partition $\lambda$
counted by $B(x_1,x_2,...,x_{n-1};m')$. We call the second
partition $\tau$ and color its parts with the color $2^{n-1}$.

$\bullet$ {\bf Step 1.} First, we change $\lambda$ in the
following way: for each part $k$ of $\tau$ that is less than or
equal to the number of parts of $\lambda$, add $1$ to each of the
first $k$ parts of $\lambda$ and then add $2^{n-1}$ to the
\emph{color} of the $k$th part. Here we record that we have

\begin{equation*}
\lambda_i - \lambda_{i+1} \geq
\begin{cases}
\omega(c_i) + \delta(c_i,c_{i+1}), & \text{if $c_i < 2^{n-1}$}, \\
\omega(c_i) + \delta(c_i-2^{n-1},c_{i+1}), & \text{if $c_i > 2^{n
- 1}$}.
\end{cases}
\end{equation*}
and that $\lambda_i\ge \omega(c_i)$ for $c_i\neq 2^{n-1}$. To be
concise, we shall write $\delta^*(c_i,c_{i+1})$ to mean
$\delta(c_i,c_{i+1})$ when $c_i < 2^{n-1}$ and
$\delta(c_i-2^{n-1},c_{i+1})$ when  $c_i > 2^{n-1}$.

\medskip
$\bullet$ {\bf Step 2.}  Now write the unused parts from $\tau$ in
descending order to the left of the parts from $\lambda$. Remove a
staircase, i.e., subtract $0$ from the smallest part, $1$ from the
next smallest, and so on.  Here we record that we have

\begin{equation} \label{pidifference}
\lambda_i - \lambda_{i+1} \geq
\omega(c_i) + \delta^*(c_i,c_{i+1}) - 1
\end{equation}
for $c_i\neq 2^{n-1}$.

\medskip
$\bullet$ {\bf Step 3.}  Starting from the largest part $k$ with
color $2^{n-1}$, we insert $k$ into the partition $\lambda$ as the
part $\lambda_i$ so that $\lambda_i - \lambda_{i+1} \geq 0$ with
$i$ minimal.  Since $\omega(c_i) = \omega(2^{n-1})=1$ and
$\delta(2^{n-1},c_{i+1})=0$, this condition is the same as
$\lambda_i-\lambda_{i+1}\ge \omega(c_i) + \delta(c_i,c_{i+1}) -
1$.  The minimality of $i$ guarantees that $\lambda_{i-1} -
\lambda_i \geq 1$.  Hence it is possible that
\begin{equation} \label{bad}
\lambda_{i-1} -\lambda_i < \omega(c_{i-1}) + \delta(c_{i-1},c_i) -
1,
\end{equation}
and in this case we shall execute a redistribution of colors
between $\lambda_i$ and $\lambda_{i-1}$.  Specifically, if
$\lambda_{i-1} - \lambda_i = j$, where $j \geq 1$ and
$\omega(c_{i-1})>1+j-\delta(c_{i-1},c_i)$, we form two new parts
$\tilde{\lambda}_{i-1}$ and $\tilde{\lambda}_i$ with colors
$\tilde{c}_{i-1}$ and $\tilde{c}_{i}$ by taking the first $j$
smallest powers of two from the color $c_{i-1}$, adding them to
the color $2^{n-1}$ to get the color $\tilde{c}_{i-1}$, and
letting $\tilde{c}_{i}$ be what is left of $c_{i-1}$.

Some comments on this change of colors are in order.  First, note
that $\tilde{c}_i \neq 2^{n-1}$, $z(\tilde{c}_i)=z(c_{i-1})$,
$v(\tilde{c}_{i-1})=v(c_{i-1})$, and
$\delta(\tilde{c}_{i-1}-2^{n-1},\tilde{c}_{i}) = 1$.  Second,
\begin{eqnarray*}
\tilde{\lambda}_{i-1} - \tilde{\lambda}_{i} &=& j \\
&=& \omega(\tilde{c}_{i-1}) +
\delta(\tilde{c}_{i-1},\tilde{c}_{i}) - 1,
\end{eqnarray*}
since the fact that $\tilde{c}_{i-1} > 2^{n-1}$ implies that
$\delta(\tilde{c}_{i-1},\tilde{c}_{i}) = 0$.  However, as
$\delta(\tilde{c}_{i-1} - 2^{n-1},\tilde{c}_{i}) = 1$, this
difference $j$ between $\tilde{\lambda}_{i-1}$ and
$\tilde{\lambda}_i$ is one less than the minimum difference
guaranteed by the second case of \eqref{pidifference}, which makes
the change of colors bijective - one can always identify when it
has taken place.  We are also guaranteed that any further
occurrences of $k$ in color $2^{n-1}$ may now be inserted without
any problem.

Next, since we are in the case of (\ref{bad}) when adding in the
part $k$ as $\lambda_i$, we may observe that $c_{i-1}$ is not
$2^{n-1}$. Moreover $c_{i+1}=2^{n-1}$ would contradict the fact
that we start with the largest part $k$. So we had $\lambda_{i-1}
- \lambda_{i+1} \geq \omega(c_{i-1}) + \delta^*(c_{i-1},c_{i+1}) -
1$ according to \eqref{pidifference}. Hence we may deduce that
\begin{eqnarray*}
\tilde{\lambda}_{i} - \lambda_{i+1}
&=& \tilde{\lambda}_i - \lambda_{i-1} + \lambda_{i-1} - \lambda_{i+1} \\
&\geq& \omega(c_{i-1}) - j + \delta^*(c_{i-1},c_{i+1}) - 1 \\
&=& \omega(\tilde{c}_i) + \delta^*(\tilde{c}_i,c_{i+1}) - 1,
\end{eqnarray*}
as $z(\tilde{c}_i)=z(c_{i-1})$ and $\tilde{c}_i\neq 2^{n-1}$.

Note also that there is no change in the required difference
between $\lambda_{i-2}$ and $\tilde{\lambda}_{i-1}$ because
$\tilde{\lambda}_{i-1} = \lambda_{i-1}$ and $v(\tilde{c}_{i-1}) =
v(c_{i-1})$ implies that $\delta({c}_{i-2},c_{i-1})=
\delta({c}_{i-2},\tilde{c}_{i-1})$.   \\

We continue this procedure until all the parts of color $2^{n-1}$
are inserted in $\lambda$.

\medskip
$\bullet$ {\bf Step 4.} Now all the required differences are

\begin{equation} \label{lambdadifffinal}
\lambda_i - \lambda_{i+1} \geq \omega(c_i) + \delta(c_i,c_{i+1}) -
1,
\end{equation}
and we can add back the staircase $0,1,2,...$ so that these
difference conditions become those of condition $(iii)$ in the
theorem. Conditions $(i)$ and $(ii)$ are straightforward. This
establishes that $A(x_1,\ldots,x_n;m) = B(x_1,\ldots,x_n;m)$.

\qed

\noindent {\bf Remark.} It will prove useful to take advantage of
the symmetry in the partitions counted by $A(x_1,\ldots,x_n;m)$
and $B(x_1,\ldots,x_n;m)$ to slightly extend Theorem
\ref{lovejoy1}. Given a permutation $\sigma \in S_n$, take a
partition $\lambda$ counted by $B(x_1,\ldots,x_n;m)$ and create a
new partition $\tilde{\lambda}$ by setting
$\tilde{\lambda}_i=\lambda_i$ and changing $c_i$ to $\sigma(c_i)$.
Here $\sigma(c_i)$ is defined in the obvious way, by permuting the
powers of $2$ occurring in the binary representation of $c_i$. Now
$\omega(c_i)$ hasn't changed so the new difference condition on
$\tilde{\lambda}$ is
\begin{equation} \label{deltasigma}
\tilde\lambda_i - \tilde\lambda_{i+1} \geq \omega(\tilde{c}_i) +
\delta(\sigma^{-1}(\tilde{c}_i),\sigma^{-1}(\tilde{c}_{i+1})).
\end{equation}
The mapping $\lambda \to \tilde{\lambda}$ is easily reversible so
we have $B(x_1,\ldots,x_n;m) =
B_{\sigma}(x_{\sigma(1)},...,x_{\sigma(n)};m)$, where
$B_{\sigma}(x_1,\ldots,x_n;m)$ denotes the number of partitions
$\lambda$ of $m$ into parts that come in $2^{n}-1$ colors and
satisfy conditions $(i)$ and $(ii)$ of Theorem \ref{lovejoy1} as
well as \eqref{deltasigma}. From the definition of $A(x_1,\ldots
,x_n;m)$, it is obvious that for any permutation
$\tau=(\tau(1),\ldots ,\tau(n))$ in $S_n$,
$$
A(x_1,\ldots ,x_n;m)=A(x_{\tau(1)},\ldots ,x_{\tau(n)};m).
$$
Hence we have
\begin{corollary} \label{lovejoy1.5}
$B_{\sigma}(x_{\tau(1)},...,x_{\tau(n)};m) = A(x_1,\ldots,x_n;m)$.
\end{corollary}

To conclude this section we provide yet another example of the
proof that $A(x_1,\ldots,x_n;m) = B(x_1,\ldots,x_n;m)$, this time
in the case $n = 4$. We start with a partition $\lambda =
(16_3,14_7,11_5,8_1,6_2,5_7,1_1)$ counted by $B(6,4,3;61)$ and a
partition $\tau = (22_8,19_8,18_8,11_8,7_8,4_8,2_8,1_8)$ of $84$
into $8$ distinct parts of color $8$.  Then we follow the steps of
the bijection:
\begin{eqnarray*}
(\tau,\lambda) &\iff&
((22_8,19_8,18_8,11_8),(20_{11},17_{15},13_5,10_9,7_{2},6_{7},2_9))
\hskip.5 in
\text{(Step 1)}\\
&\iff&
((12_8,10_8,10_8,4_8),(14_{11},12_{15},9_5,7_9,5_{2},5_{7},2_9))
\hskip.5 in
\text{(Step 2)}\\
&\iff&
((10_8,10_8,4_8),(14_{11},12_{8},12_{15},9_5,7_9,5_{2},5_{7},2_9))
\hskip.5 in
\text{(Step 3)}\\
&\iff& ((10_8,4_8),(14_{11},12_{8},\underbrace{12_{15},10_8},9_5,7_9,5_{2},5_{7},2_9)) \\
&\iff& ((10_8,4_8),(14_{11},12_{8},12_{11},10_{12},9_5,7_9,5_{2},5_{7},2_9)) \\
&\iff& ((4_8),(14_{11},12_{8},12_{11},10_8,10_{12},9_5,7_9,5_{2},5_{7},2_9)) \\
&\iff& (\epsilon,(14_{11},12_{8},12_{11},10_8,10_{12},9_5,7_9,5_{2},\underbrace{5_{7},4_8},2_9)) \\
&\iff& (\epsilon,(14_{11},12_{8},12_{11},10_8,10_{12},9_5,7_9,5_{2},{5_{9},4_6},2_9)) \\
&\iff&
(24_{11},21_{8},20_{11},17_8,16_{12},14_5,11_9,8_{2},{7_{9},5_6},2_9)
\hskip.5 in \text{(Step 4)}
\end{eqnarray*}
Notice that the result is a partition counted by $B(6,4,3,8;145)$,
as expected.

\section{Proofs of Theorems \ref{lovejoy2} and \ref{lovejoy3}}

The proofs of these theorems follow exactly the same steps as the
proof of Theorem \ref{lovejoy1}, but here we will pay special
attention to the behavior of the parts with odd color.  We begin
with Theorem \ref{lovejoy2}. When $n=1$, we just have a partition
$\mu_1$ into parts having color $1$ that are congruent to $r$
modulo $M$, and the conditions $(i)$ and $(ii)$ of the theorem are
trivial.  Now suppose that the theorem is true for $n-1$, let
$\lambda$ be a partition counted by $B_{r,M}(x_1,...,x_{n-1};m')$,
and let $\tau$ be a partition of $m-m'$ into distinct parts, all
with the color $2^{n-1}$. As we apply the bijection, the part
$\lambda_{S}$ will increase by $1$ each time its color or the
color of a smaller part increases by $2^{n-1}$ in Step 1. It will
also ultimately increase by $1$ in Step 4 if a part of color
$2^{n-1}$ is inserted as a part $\lambda_{S + k}$ in Step 3. But
in these cases $\sum_{\ell = S}^{s} \omega_e(c_{\ell})$ also
increases by $1$.  If a redistribution of colors should take place
between $\lambda_{S}$ and $\lambda_{S + 1}$, then
$\tilde{\lambda}_{S}$ remains the smallest part with odd color,
and $\sum_{\ell = S}^{s} \omega_e(c_{\ell})$ does not change.
Hence we have condition $(i)$ of Theorem \ref{lovejoy2}.

For condition $(ii)$, the difference between two parts $\lambda_i$
and $\lambda_j$ with odd color will increase by $1$ every time in
Step 1 that the color of $\lambda_k$ increases by $2^{n-1}$ for $i
\leq k < j$.  This difference will also ultimately increase by $1$
in Step 4 for each time that a part with color $2^{n-1}$ is
inserted between $\lambda_i$ and $\lambda_j$ in step 3.  But in
these cases $\sum_{\ell = i}^{j-1} \omega_e(\lambda_{\ell})$ will
also increase by $1$. If any part $\lambda_i$ with odd color is
affected by a rearrangement of colors, $\tilde{\lambda}_i$ remains
odd.  If $\lambda_j$ has odd color with $j > i$, then $\sum_{\ell
= i}^{j-1} \omega_e(\lambda_{\ell})$ is not affected by this
rearrangement, and if $j < i$, neither is $\sum_{\ell = j}^{i-1}
\omega_e(\lambda_{\ell})$.  This guarantees condition $(ii)$ and
completes the proof of this part of the theorem.

\qed

We turn to Theorem \ref{lovejoy3}.  When $n=1$, we just have a
partition $\mu_1$ into parts having color $1$ that differ by at
least $M$, and the extra condition of the theorem is trivial.  Now
suppose that the theorem is true for $n-1$, let $\lambda$ be a
partition counted by $B_M(x_1,...,x_{n-1};m')$, and let $\tau$ be
a partition of $m - m'$ into distinct parts, all with the color
$2^{n-1}$. As we apply the steps of the bijection, the difference
between any two parts with odd color $\lambda_i$ and $\lambda_j$
increases by $1$ for every time in Step 1 that the color of
$\lambda_k$ increases by $2^{n-1}$ for $i \leq k < j$. This
difference will also ultimately increase by $1$ in Step 4 for each
time that a part with color $2^{n-1}$ is inserted between
$\lambda_i$ and $\lambda_j$ in step 3.  But these are precisely
the cases where $\sum_{\ell = i}^{j-1} \omega_e(c_{\ell})$
increases by $1$.  If any part $\lambda_i$ with odd color is
affected by a rearrangement of colors, $\tilde{\lambda}_i$ remains
odd.  If $\lambda_j$ has odd color with $j < i$, then $\lambda_j -
\tilde{\lambda}_i$ is not affected by this rearrangement, and
neither is $\sum_{\ell = j}^{i-1} \omega_e(c_{\ell})$.  If $j
> i$, then both $\tilde{\lambda}_i - \lambda_j$ and $\sum_{\ell =
i}^{j-1} \omega_e(c_{\ell})$ are unchanged. This guarantees the
extra condition and completes the proof of this part of the
theorem.

\qed

\section{Partition Identities} \setcounter{equation}{0} \setcounter{theorem}{0}

We begin our discussion of partition identities by proving Theorem
\ref{lovejoy1conj}.

\noindent {\bf Proof of Theorem \ref{lovejoy1conj}.} For each part
$\lambda_i$ of a partition counted by $B(x_1,\ldots,x_n;m)$ with
color $c_i = 2^{j_1} + \cdots +2^{j_k}$ with $j_1<\cdots <j_k$, we
draw its Ferrers diagram, that is we write a row of $\lambda_i$
boxes  and we add subscripts on the last $k$ boxes. Specifically,
the last box has subscript $j_1+1$, the next to last has $j_2+1$,
and so on. Conjugating this diagram and interpreting a subscript
as the color of a column gives a
partition counted by $C(x_1,\ldots,x_n;m)$. \qed\\

\noindent{\bf Example.} Let us take $n=3$ and
$\lambda=(10_7,6_5,4_1,1_2)$ in $B(3,2,2;21)$. The Ferrers diagram
with the subscripts is
$$
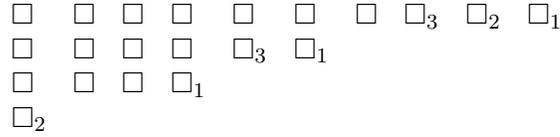
\begin{array}{llllllllll}
\Box & \Box & \Box & \Box & \Box & \Box & \Box & \Box_3 & \Box_2 & \Box_1 \\
\Box & \Box & \Box & \Box & \Box_3 & \Box_1&&&&\\
\Box & \Box & \Box & \Box_1&&&&&&\\
\Box_2 &&&&&&&&&\\
\end{array}
$$
Now we read the columns and get
$(4_2,3,3,3_1,2_3,2_1,1,1_3,1_2,1_1)$ which
is in $C(3,2,2;21)$.\\

We note, as was done in \cite[p. 25]{Al1}, that making the
substitutions $q \to q^N$ and $y_j \to q^{a_{j} - N}$ in
\eqref{product} and applying Theorem \ref{lovejoy1conj} gives the
Andrews-Olsson identities referred to in the introduction:

\begin{theorem} [Andrews-Olsson] \label{andrewsolsson}
Let $N$ be a positive integer and let $A =\{a_1,a_2,\dots,a_n\}$
be a set of distinct positive integers arranged in increasing
order with $a_n < N$.  Let $P_1(A;N;m)$ denote the number of
partitions of $m$ into distinct parts each congruent to some $a_i
\pmod{N}$.  Let $P_2(A;N;m)$ denote the number of partitions of
$m$ into parts $\equiv 0$ or some $a_i$ modulo $N$ such that only
parts divisible by $N$ may repeat, the smallest part is less than
$N$, and the difference between parts is $\leq N$, with strict
inequality if either part is divisible by $N$.  Then $P_1(A;N;m) =
P_2(A;N;m)$.
\end{theorem}


We also note, before continuing, that there are conjugate versions
of Theorems \ref{lovejoy2} and \ref{lovejoy3} as well, the extra
conditions on the differences between parts of odd color
translating under conjugation to conditions on the number of parts
occurring between two parts of color $1$.

Next we discuss how a theorem of Bressoud \cite{Br.5}, which
generalizes some results of G\"ollnitz \cite{Go1}, is contained in
Theorem \ref{lovejoy2}.

\begin{theorem} [Bressoud]
Given positive integers $n,k$, and $r$ satisfying $1 \leq r < 2k$
and $r \neq k$, let $G_{r,k}(n)$ denote the number of partitions
of $n$ into distinct parts congruent to $r$, $k$, or $2k$ modulo
$2k$, and let $H_{r,k}(n)$ denote the number of partitions of $n$
into parts congruent to $r$ or $k$ modulo $k$ with minimal
difference $k$, minimal difference $2k$ between parts congruent to
$r$ modulo $k$, and, if $r > k$, with the smallest part greater
than or equal to $k$.  Then $G_{r,k}(n) = H_{r,k}(n)$.
\end{theorem}

\noindent {\bf Proof.} This theorem is a special case of Theorem
\ref{lovejoy2} when $n,M = 2$.  We supply the details for $r < k$,
where we take $R = 1$ in Theorem \ref{lovejoy2}; for the other
case, which is similar, use $R = 2$. To begin, let $\lambda$ be a
partition counted by $B_{1,2}(x_1,x_2;m)$.  The important
observation is that the extra conditions $(i)$ and $(ii)$ in
Theorem \ref{lovejoy2} ensure that we can ``drop" the color $2$
from the subscripts without losing any information. The color $1$
remains $1$, the color $2$ is dropped, and the color $3$ becomes
$1$. This operation corresponds to setting $a_2 = 1$ in the
product $(-a_1q;q^2)_{\infty}(-a_2q;q)_{\infty}.$  For example,
the partition $(13_{3},10_1,8_2,7_1,5_3,3_2)$ becomes
$(13_1,10_1,8,7_1,5_1,3)$. One verifies that the difference
conditions on parts of $\lambda$ become $\lambda_i - \lambda_{i+1}
\geq 2$, if $c_i = 1$, and $\lambda_i - \lambda_{i+1} \geq 1$ if
$c_i$ is uncolored.  To finish, we replace $q$ by $q^k$ and $a_1$
by $q^{r-k}$ in the product $(-a_1q;q^2)_{\infty}(-q;q)_{\infty}$,
which corresponds to replacing a part $j_1$ by $(j-1)k + r$ and a
part $j$ by $kj$ in $\lambda$.  The difference conditions are now
precisely those of Bressoud's theorem.    \qed

In the late 1960's Andrews \cite{An1,An2} proved two superficially
similar generalizations of Schur's theorem, which we need some
notation to state. Consider again a set $A = \{a_1,a_2,\cdots,a_n
\}$ of $n$ distinct positive integers, this time satisfying
$\sum_{i=1}^{k-1} a_i < a_k$ for all $1 \leq k \leq n$. Fix an
integer $N$ such that $N \ge \sum_{i=1}^n a_i$. Denote the set of
$2^n - 1$ (necessarily distinct) possible sums of distinct
elements of $A$ by $A'$ and its elements by $\alpha_1 < \alpha_2 <
\cdots < \alpha_{2^n-1}$. For any natural number $N \geq
\alpha_{2^n-1}$ define $A_N$ (resp. $A'_N, -A_N, -A'_N$) to be the
set of all natural numbers congruent to some $a_i$ (resp.
$\alpha_i, -a_i, -\alpha_i$) modulo $N$. For $x \in A'$ let
$\omega_A(x)$ denote the number of terms in the defining sum of
$x$ and let $v_A(x)$ (resp. $z_A(x)$) be the smallest (resp.
largest) $a_i$ appearing in this sum. Moreover, for $x,y\in A'$,
we define $\delta_A(x,y)=1$ if $z_A(x)< v_A(y)$ and 0 otherwise.
Finally, let $\beta_N(\ell)$ be the least positive residue of
$\ell$ modulo $N$.

\begin{theorem} [Andrews, \cite{An2,An3}] \label{andrews1}
Let $D(A_N;x_1,\ldots,x_n;m)$ denote the number of partitions of
$m$ into distinct parts taken from $A_N$ where there are $x_r$
parts equivalent to $a_r$ modulo $N$. Let
\newline $E(A'_N;x_1,\ldots,x_n;m)$ denote the number of partitions of
$m$ into parts taken from $A'_N$ of the form $\lambda_1 + \cdots +
\lambda_s$ such that  (i) $x_r$ is the number of $\lambda_i$ such
that $\beta_N(\lambda_i)$ uses $a_r$ in its defining sum, and
$(ii)$
$$\lambda_i - \lambda_{i+1} \geq N \omega_A(\beta_N(\lambda_{i+1}))
+ v_A(\beta_N(\lambda_{i+1})) - \beta_N(\lambda_{i+1}).$$  Then
$D(A_N;x_1,\ldots,x_n;m) = E(A'_N;x_1,\ldots,x_n;m)$.
\end{theorem}

\begin{theorem} [Andrews, \cite{An1,An3}] \label{andrews2}
Let $F(-A_N;x_1,\ldots,x_n;m)$ denote the number of partitions of
$m$ into distinct parts taken from $-A_N$ where there are $x_r$
parts equivalent to $-a_r$ modulo $N$. Let
$G(-A'_N;x_1,\ldots,x_n;m)$ denote the number of partitions of $m$
into parts taken from $-A'_N$ of the form $\lambda_1 + \cdots
+\lambda_s$ such that $(i)$ $\lambda_i \ge
N(\omega_A(\beta_N(-\lambda_i))-1)$, (ii) $x_r$ is the number of
$\lambda_i$ such that $\beta_N(-\lambda_i)$ uses $a_r$ in its
defining sum, and (iii) $$\lambda_i - \lambda_{i+1} \geq
N\omega_A(\beta_N(-\lambda_i)) + v_A(\beta_N(-\lambda_{i})) -
\beta_N(-\lambda_{i})$$. Then $F(-A_N;x_1,\ldots,x_n;m) =
G(-A'_N;x_1,\ldots,x_n;m)$.
\end{theorem}

In the rest of the paper we first state and then prove a number of
partition theorems that extend the results of Andrews. These will
correspond to the substitutions $q \to q^N$ and $y_j \to q^{a_j -
N}$ in \eqref{product}, \eqref{product1}, and \eqref{product2},
and $q \to q^N$ and $y_j \to q^{-a_j}$ in \eqref{product}.

For the first two results, we take advantage of the symmetry in
Theorem \ref{deltasigma}.  We extend any permutation $\sigma \in
S_n$ to the integers in $A'$ in the obvious way, by letting
$\sigma(\alpha_i)=\sum_{j=1}^k a_{\sigma(i_j)}$ if
$\alpha_i=\sum_{j=1}^k a_{i_j}$.

\begin{theorem} \label{corteel1}
Let $E_\sigma(A'_N;x_1,\ldots,x_n;m)$ denote the number of
partitions of $m$ into parts taken from $A'_N$ of the form
$\lambda_1 + \cdots + \lambda_s$ such that  (i) $x_r$ is the
number of $\lambda_i$ such that $\beta_N(\lambda_i)$ uses $a_r$ in
its defining sum, and (ii)
$$
\lambda_i - \lambda_{i+1} \geq N \omega_A(\beta_N(\lambda_{i+1}))
+N\delta_A(\sigma(\beta_N(\lambda_{i})),\sigma(\beta_N(\lambda_{i+1})))
+\beta_N(\lambda_{i})- \beta_N(\lambda_{i+1}).$$ Then
$D(A_N;x_1,\ldots,x_n;m)=E_\sigma(A'_N;x_1,\ldots,x_n;m)$.
\end{theorem}

\begin{theorem} \label{corteel2}
Let $G_\sigma(-A'_N;x_1,\ldots,x_n;m)$ denote the number of
partitions of $m$ into parts taken from $-A'_N$ of the form
$\lambda_1 + \cdots +\lambda_s$ such that $(i)$ $\lambda_i \ge
N(\omega_A(\beta_N(-\lambda_i))-1)$, (ii) $x_r$ is the number of
$\lambda_i$ such that $\beta_N(-\lambda_i)$ uses $a_r$ in its
defining sum, and (iii)
$$\lambda_i - \lambda_{i+1} \geq
N\omega_A(\beta_N(-\lambda_i))
+N\delta_A(\sigma(\beta_N(-\lambda_{i})),\sigma(\beta_N(-\lambda_{i+1})))+\beta_N(-\lambda_{i+1})
- \beta_N(-\lambda_{i}).$$ Then $F(-A_N;x_1,\ldots,x_n;m) =
G_\sigma(-A'_N;x_1,\ldots,x_n;m)$.
\end{theorem}

Let us record some examples.  Suppose that $n=2$, $A=\{1,2\}$ and
$N=3$.  It is easy to see that taking $\sigma=(1,2)$ in Theorem
\ref{corteel1} or $\sigma= (2,1)$ in Theorem \ref{corteel2} gives
back (a refinement of) Schur's theorem.
On the other hand, taking $\sigma=(2,1)$ in Theorem \ref{corteel1}
or $\sigma=(1,2)$ in Theorem \ref{corteel2} gives
\begin{corollary}
The number of partitions of $m$ into distinct parts $\equiv 1,2
\pmod3$ with $x_r$ parts congruent to $r$ modulo 3 is equal to the
number of partitions $\lambda$ of $m$ such that (i) $x_r$ is the
number of parts congruent to $r$ or $3$ modulo 3 and (ii) if
$\lambda_{i}\equiv j\pmod3$ and $\lambda_{i+1}\equiv k\pmod3$ then
\[
\lambda_i-\lambda_{i+1}\ge \left\{ \begin{array}{ll} 3,& {\rm if}\ k=1\ {\rm and} \ j=1,3,\\
 7, & {\rm if}\ k=1\ {\rm and} \ j=2,\\
 2, & {\rm if}\ k=2,\\
 4, & {\rm if}\  k=3.\\
\end{array}\right.
\]
\end{corollary}
\noindent This is equivalent to the case $S_1 = S_4$ of
\cite[Theorem 8]{Al-Go2}.

For a more complicated example, let $n=3$, $A=\{1,2,4\}$ and
$N=7$. Taking $\sigma=(3,2,1)$ in Theorem \ref{corteel1} gives
\begin{corollary}
The number of partitions of $m$ into distinct parts $\equiv 1,2,4
\pmod7$ with $x_r$ parts congruent to $2^{r-1}$ modulo 7  is equal
to the number of partitions $\lambda$ of $m$ such that $(i)$ $x_r$
is the number of $\lambda_i$ such that $\lambda_i\pmod7$ uses
$2^{r-1}$ in its defining sum, and (ii) if $\lambda_{i}\equiv
j\pmod7$ and $\lambda_{i+1}\equiv k\pmod7$,
\[
\lambda_i-\lambda_{i+1}\ge\left\{\begin{array}{ll}
 7, & {\rm if}\ k=1\ {\rm and}\ j=1,3,5,7,\\
 13, & {\rm if}\ k=1\ {\rm and}\ j=2,4,6,\\
 6, & {\rm if}\ k=2\ {\rm and}\ j\neq 4,\\
 16, & {\rm if}\ k=2\ {\rm and}\ j=4,\\
 12, & {\rm if}\ k=3\ {\rm and}\ j\neq 4,\\
 22, & {\rm if}\ k=3\ {\rm and}\ j=4,\\
 4, & {\rm if}\ k=4,\\
 10, & {\rm if}\ k=5,\\
 9, & {\rm if}\ k=6,\\
 15, & {\rm if}\  k=7.\\
\end{array}\right.
\]
\end{corollary}

Although it may not yet be clear, we shall see that Theorem
\ref{andrews1} is the case $\sigma = Id$ of Theorem \ref{corteel1}
and Theorem \ref{andrews2} is the case $\sigma = (n,n-1,...,1)$ of
Theorem \ref{corteel2}.   While these theorems take advantage of
the symmetry in Theorem \ref{lovejoy1.5}, this symmetry does not
persist in Theorems \ref{lovejoy2} and \ref{lovejoy3}, from which
our last two theorems follow. For $\ell \in A'$ let
$\omega_{A,1}(\ell)$ denote the number of terms not
equal to $a_1$ in the defining sum of $\ell$.\\

\begin{theorem} \label{corteel3}
Let $D_{R,M}(A_{N};x_1,\ldots,x_n;m)$ denote the number of
partitions counted by \newline $D(A_{N};x_1,\ldots,x_n;m)$ such
that the parts that are $a_1$ modulo $N$ are $((R-1)N+a_1)$ modulo
$MN$. Let $E_{R,M}(A'_{N};x_1,\ldots,x_n;m)$ denote the number of
partitions of $m$ counted by $E(A'_{N};x_1,\ldots,x_n;m)$ of the
form $\lambda_1 + \cdots + \lambda_s$ such that (i) if $\lambda_S$
is the smallest part such that $\beta_N(\lambda_S)$ uses $a_1$ in
its defining sum, then
$$
\lambda_S\equiv
N(R-\omega_A(\beta_N(\lambda_S)))+\beta_N(\lambda_S)+
N\sum_{\ell=S}^{s}\omega_{A,1}(\beta_N(\lambda_{\ell})) \pmod{MN},
$$ and (ii) if $i<j$ and $\beta_N(\lambda_i)$ and
$\beta_N(\lambda_j)$ use $a_1$ in their defining sums, then
$$
\lambda_i-\lambda_j\equiv
N(-\omega_A(\beta_N(\lambda_i))+\omega_A(\beta_N(\lambda_{j})))+
\beta_N(\lambda_i)-\beta_N(\lambda_{j})
+N\sum_{\ell=i}^{j-1}\omega_{A,1}(\beta_N(\lambda_{\ell}))
\pmod{MN}.
$$
Then $D_{R,M}(A_{N};x_1,\ldots,x_n;m) =
E_{R,M}(A'_{N};x_1,\ldots,x_n;m)$
\end{theorem}

For $R=1$, $M=2$, $N=3$, and $A = \{1,2\}$, this translates to
\begin{corollary}
Let $D_{1,2}(A_3;x_1,x_2;m)$ denote the number of partitions of
$m$ counted by \newline $D(A_3;x_1,x_2;m)$ where the parts that
are 1 modulo 3 are $1$ modulo $6$, and let
$E_{1,2}(A'_3;x_1,x_2;m)$ denote the number of partitions of $m$
counted by $E(A'_3;x_1,x_2;m)$ of the form $\lambda_1 + \cdots +
\lambda_s$ such that
if $\lambda_S$ is the smallest part that is $\equiv 1$ or
$3\pmod3$, then
$(i)$ if $\lambda_S\equiv 1\pmod3$ then $ \lambda_S\equiv
1+3\sum_{\ell=S}^{s}\omega_{A,1}(\beta_3(\lambda_{\ell}))
\pmod{6}$, 
$(ii)$ if $\lambda_S\equiv 3\pmod3$ then $ \lambda_S\equiv
3\sum_{\ell=S}^{s}\omega_{A,1}(\beta_3(\lambda_{\ell})) \pmod{6},
$
and $(iii)$ if $i<j$ and $\lambda_i,\lambda_j\equiv 1,3\pmod3$,
then
$$\lambda_i-\lambda_j\equiv
\frac{\beta_3(\lambda_i)-\beta_3(\lambda_j)}{2}+
3\sum_{\ell=i}^{j-1}\omega_{A,1}(\beta_3(\lambda_{\ell}))\pmod6.
$$
Then $D_{1,2}(A_3;x_1,x_2;m)=E_{1,2}(A'_3;x_1,x_2;m)$.
\label{beforeagain}
\end{corollary}

\noindent{\bf Example.}  For $m=18$, we have $D(A_3;1,1;18)=6$, the
relevant partitions being $(17,1)$, $(14,4)$, $(11,7)$, $(10,8)$,
$(13,5)$ and $(16,2)$. But only three of them satisfy the conditions
of the previous corollary, i.e., that the parts that are 1 modulo 3
are 1 modulo 6. These partitions are $(17,1)$, $(11,7)$ and
$(13,5)$. Therefore $D_{1,2}(A_3;1,1;18)=3$. Now $E(A'_3;1,1;18)=6$,
the relevant partitions being $(18)$, $(17,1)$, $(16,2)$, $(14,4)$,
$(13,5)$ and $(11,7)$. But three of them violate the condition on
$\lambda_S$, namely $(18)$ violates $(ii)$, and $(14,4)$ and
$(13,5)$ violate $(i)$.
So, $E_{1,2}(A'_{3},1,1;18)=3=D_{1,2}(A_3;1,1;18)$.\\

\begin{theorem} \label{corteel4} Let $D_{M}(A_{N};x_1,\ldots,x_n;m)$
denote the number of partitions counted by
\newline $D(A_{N};x_1,\ldots,x_n;m)$ such that the parts equivalent to
$a_1$ modulo $N$ differ at least by $MN$. Let
$E_{M}(A'_{N};x_1,\ldots,x_n;m)$ denote the number of partitions
of $m$ counted by $E(A'_{N};x_1,\ldots,x_n;m)$ of the form
$\lambda_1 + \cdots + \lambda_s$ such that  if $i<j$ and
$\beta_N(\lambda_i)$ and $\beta_N(\lambda_j)$ use $a_1$ in their
defining sums, then
$$
\lambda_i-\lambda_j\ge
MN+N(-\omega_A(\beta_N(\lambda_i))+\omega_A(\beta_N(\lambda_{j})))
+\beta_N(\lambda_i)-\beta_N(\lambda_{j})
+N\sum_{\ell=i}^{j-1}\omega_{A,1}(\beta_N(\lambda_{\ell})).
$$
Then $D_{M}(A_{N};x_1,\ldots,x_n;m) =
E_{M}(A'_{N};x_1,\ldots,x_n;m)$
\end{theorem}

Let $M=2$, $N=3$, and $A = \{1,2\}$~:
\begin{corollary}
Let $D_{2}(A_3;x_1,x_2;m)$ denote the number of partitions of $m$
counted by \newline $D(A_3;x_1,x_2;m)$ such that the parts
equivalent to $1$ modulo $3$ differ by at least $6$ and let
\newline $E_{2}(A'_3;x_1,x_2;m)$ denote the number of partitions of $m$
counted by $E(A'_3;x_1,x_2;m)$
such that if $i<j$ and $\lambda_i,\lambda_j\equiv 1,3\pmod3$ then
$\lambda_i-\lambda_j\ge
5+3\sum_{\ell=i}^{j-1}\omega_{A,1}(\beta_3(\lambda_{\ell}))$.
Then $D_{2}(A_3;x_1,x_2;m)=E_{2}(A'_3;x_1,x_2;m)$.
\label{beforeagain2}
\end{corollary}
\noindent{\bf Example.} For $m=22$, we have
$D_{2}(A_{3};2,1;22)=8$, with the relevant partitions being
$(16,4,2)$, $(16,5,1)$,  $(13,7,2)$, $(13,5,4)$ $(13,8,1)$,
$(14,7,1)$, $(10,8,4)$  and $(11,10,1)$. Also
$E_{2}(A'_{3},2,1;22)=8$, with the relevant partitions being
$(18,4)$, $(19,3)$, $(15,7)$, $(16,6)$, $(16,5,1)$, $(14,7,1)$,
$(13,7,2)$ and $(13,8,1)$. \\

We now turn to the proofs of all of these theorems.

\noindent{\bf Proof of Theorem \ref{corteel1}.} Fix a permutation
$\sigma \in S_n$. In a partition counted by $A(x_1,\ldots,x_n;m)$,
we replace a part of size $k$ from $\mu_j$ by $N(k - 1) + a_j$.
Then each $\lambda_i$ in the corresponding partition $\lambda$
counted by $B_{\sigma^{-1}}(x_1,\ldots,x_n;m)$ is replaced by
$\tilde\lambda_i=N(\lambda_i - \omega_A(\alpha_{c_i})) +
\alpha_{c_i}$. This corresponds to replacing $q$ by $q^N$ and
$y_j$ by $q^{a_{j} - N}$ in \eqref{product}. Suppose that before
this replacement, we had $\lambda_{i+1} = \lambda_i - \omega(c_i)
- \delta(\sigma(c_i),\sigma(c_{i+1})) - d$, where $d \geq 0$. Then
after the replacement we have
$$
\tilde{\lambda}_i = N(\lambda_i - \omega_A(\alpha_{c_i})) +
\alpha_{c_i}
$$
and
$$
\tilde{\lambda}_{i+1} = N(\lambda_i - \omega_A(\alpha_{c_i}) -
\delta(\sigma(c_i),\sigma(c_{i+1})) - d -
\omega_A(\alpha_{c_{i+1}})) + \alpha_{c_{i+1}}.
$$
Since
$\delta(\sigma(c_i),\sigma(c_{i+1}))=\delta_A(\alpha_{\sigma(c_i)},\alpha_{\sigma(c_{i+1})})$,
we get
$$
\tilde{\lambda}_i - \tilde{\lambda}_{i+1} = N(d +
\delta_A(\alpha_{\sigma(c_i)},\alpha_{\sigma(c_{i+1}))} +
\omega_A(\alpha_{c_{i+1}})) + \alpha_{c_i} - \alpha_{c_{i+1}}.
$$
Note that $\beta_N(\tilde{\lambda_i})=\alpha_{c_i}$ and that
$\sigma(\beta_N(\tilde{\lambda_i}))=\alpha_{\sigma(c_i)}$ for all
$i$ and the
theorem follows. \qed \\

\noindent{\bf Proof of Theorem \ref{corteel2}.} In a partition
counted by $A(x_1,\ldots,x_n;m)$, we replace a part of size $k$
from $\mu_j$ by $Nk-a_j$.  Each $\lambda_i$ in the corresponding
partition $\lambda$ counted by $B_{\sigma^{-1}}(x_1,\ldots,x_n;m)$
is replaced by $\tilde\lambda=N\lambda_i - \alpha_{c_i}$. This
corresponds to replacing $q$ by $q^N$ and $y_j$ by $q^{-a_{j}}$ in
\eqref{product}. Suppose that before this replacement, we had
$\lambda_{i+1} = \lambda_i - \omega(c_i) -
\delta(\sigma(c_i),\sigma(c_{i+1})) - d$, where $d \geq 0$. Then
after the replacement we have
$$
\tilde{\lambda}_i = N\lambda_i - \alpha_{c_i}
$$
and
$$
\tilde{\lambda}_{i+1} = N(\lambda_i - \omega_A(\alpha_{c_i}) -
\delta(\sigma(c_i),\sigma(c_{i+1})) - d) - \alpha_{c_{i+1}}.
$$
Since
$\delta(\sigma(c_i),\sigma(c_{i+1}))=\delta_A(\alpha_{\sigma(c_i)},\alpha_{\sigma(c_{i+1})})$,
we get
$$
\tilde{\lambda}_i - \tilde{\lambda}_{i+1} = N(d +
\delta_A(\alpha_{\sigma(c_i)},\alpha_{\sigma(c_{i+1})}) +
\omega_A(\alpha_{c_{i}})) + \alpha_{c_i} - \alpha_{c_{i+1}}.
$$
Note that $\beta_N(-\tilde{\lambda_i})=\alpha_{c_i}$ and that
$\sigma(\beta_N(-\tilde{\lambda_i}))=\alpha_{\sigma(c_i)}$ for all
$i$ and the
theorem follows. \qed\\

For general $n$, it may not be clear that Theorems \ref{corteel1}
and \ref{corteel2} are indeed generalizations of Theorems
\ref{andrews1} and \ref{andrews2}. But it becomes clear with the
following lemma.

\begin{lemma}
For $x,y\in A'$,
\[
N+v_A(y)>N\delta_A(x,y)+x \ge v_A(y)
\]
\label{key}
\end{lemma}
\noindent{\bf Proof.} We consider two cases.  First, if
$\delta_A(x,y) = 0$, then $z_A(x) \geq v_A(y)$. The first
inequality is trivial as $N>x$. The second follows from the fact
that $x \geq z_A(x) \geq v_A(y)$. On the other hand, if
$\delta(x,y) = 1$, then $z_A(x)<v_A(y)$. The second inequality is
trivial as $N>v_A(y)$. The first follows from the fact that
$v_A(y)>x$ if $v_A(y) > z_A(x)$.    \qed\\

Now we can prove

\begin{corollary}
Theorem \ref{andrews1} is Theorem \ref{corteel1} with
$\sigma=(1,2,\ldots,n)$.
\end{corollary}
\noindent{\bf Proof.} As $\sigma=(1,2,\ldots,n)$, we have
$\delta_A(\sigma(\beta_N(\lambda_{i})),\sigma(\beta_N(\lambda_{i+1})))=
\delta_A(\beta_N(\lambda_{i})),\beta_N(\lambda_{i+1})$. Lemma
\ref{key} gives
\begin{equation*}
N + v_A(\beta_N(\lambda_{i+1})) >
N\delta_A(\beta_N(\lambda_{i}),\beta_N(\lambda_{i+1})) +
\beta_N(\lambda_{i}) \geq v_A(\beta_N(\lambda_{i+1})).
\end{equation*}
This shows that the minimal difference in Theorem \ref{corteel1}
is at least the one claimed in Theorem \ref{andrews1}, but not
greater.\qed

\begin{corollary}
Theorem \ref{andrews2} is Theorem \ref{corteel2} with
$\sigma=(n,n-1,\ldots,1)$.
\end{corollary}
\noindent{\bf Proof.} As $\sigma=(n,n-1,\ldots,1)$, we have
$\delta_A(\sigma(\beta_N(-\lambda_{i})),\sigma(\beta_N(-\lambda_{i+1})))=
\delta_A(\beta_N(-\lambda_{i+1}),\beta_N(-\lambda_{i}))$. Lemma
\ref{key} gives
\begin{equation*}
N + v_A(\beta_N(-\lambda_{i})) >
N\delta_A(\beta_N(-\lambda_{i+1}),\beta_N(-\lambda_{i})) +
\beta_N(-\lambda_{i+1}) \geq v_A(\beta_N(-\lambda_{i})).
\end{equation*}
This shows that the minimal difference in Theorem \ref{corteel2}
is at least the one claimed in Theorem \ref{andrews2}, but not greater.\qed\\

\noindent{\bf Proof of Theorem \ref{corteel3}.} In a partition
counted by $A_{R,M}(x_1,\ldots,x_n;m)$, we replace a part of size
$k$ from $\mu_j$ by $N(k - 1) + a_j$. In the corresponding
partition $\lambda$ counted by $B_{R,M}(x_{1},x_{2},...,x_{n};m)$,
then each $\lambda_i$ is replaced by
$\tilde{\lambda}_i=N(\lambda_i - \omega_A(\alpha_{c_i})) +
\alpha_{c_i}$. This implies that $\tilde{\lambda}\in
E(A'_{N};x_1,\ldots,x_n;m)$ as in Theorem \ref{andrews1}. Note
that $\alpha_{c_i}=\beta_N(\tilde{\lambda}_i)$,
$\omega(c_i)=\omega_A(\beta_N(\tilde\lambda_i))$ and
$\omega_e(c_{i})=\omega_{A,1}(\beta_N(\tilde\lambda_i))$. If
$\lambda_{S}$ was the smallest part with odd color, then
$\lambda_{S} \equiv R + \sum_{\ell = S}^{s} \omega_e(c_{\ell})
\pmod{M}$.  This condition is
 translated to if $\tilde\lambda_S$ is the smallest part such
that $\beta_N(\tilde\lambda_S)$  uses $a_1$ in its defining sum
$$
\tilde\lambda_{S} \equiv
RN-N\omega_A(\beta_N(\tilde\lambda_S))+\beta_N(\tilde\lambda_S)+
N\sum_{\ell = S}^{s}
\omega_{A,1}(\beta_N(\tilde\lambda_\ell))\pmod{MN}.
$$
If $\lambda_i$ and $\lambda_j$ were any two parts with odd color,
then $\lambda_i - \lambda_j \equiv \sum_{\ell = i}^{j-1}
\omega_e(\lambda_{\ell}) \pmod{M}$ and it gives $N\lambda_i -
N\lambda_j \equiv N\sum_{\ell = i}^{j-1} \omega_e(\lambda_{\ell})
\pmod{MN}$. Then $\lambda_{i}$ is changed to
$\tilde\lambda_i=N(\lambda_{i}-\omega_A(\alpha_{c_i}))+\alpha_{c_i}$
and $\lambda_{j}$ is changed to
$\tilde\lambda_j=N(\lambda_{j}-\omega_A(\alpha_{c_j}))+\alpha_{c_j}$.
We get that if $\beta_N(\tilde\lambda_i)$ and
$\beta_N(\tilde\lambda_j)$ use $a_1$ in their defining sums, then
$$
\tilde\lambda_i-\tilde\lambda_j\equiv
N(-\omega_A(\beta_N(\tilde\lambda_i))+\omega_A(\beta_N(\tilde\lambda_{j})))+
\beta_N(\tilde\lambda_i)-\beta_N(\tilde\lambda_{j})
+N\sum_{\ell=i}^{j-1}\omega_{A,1}(\beta_N(\tilde\lambda_{\ell}))
\pmod{MN}.
$$ \qed\\

\noindent{\bf Proof of Theorem \ref{corteel4}.} The proof uses the
same ideas as the proof of Theorem \ref{corteel3}. In a partition
$\lambda$ counted by $A_{M}(x_1,\ldots,x_n;m)$, we replace a part
of size $k$ from $\mu_j$ by $N(k - 1) + a_j$. In the corresponding
partition counted by $B_{M}(x_{1},x_{2},...,x_{n};m)$, then each
$\lambda_i$ is replaced by $\tilde{\lambda}_i=N(\lambda_i -
\omega_A(\alpha_{c_i})) + \alpha_{c_i}$. This corresponds to
replacing $q$ by $q^N$ and $y_j$ by $q^{a_{j} - N}$ in
\eqref{product1}. This implies that $\tilde{\lambda}\in
E(A'_{N};x_1,\ldots,x_n;m)$ as in Theorem \ref{andrews1}. Note
that $\alpha_{c_i}=\beta_N(\tilde{\lambda}_i)$,
$\omega(c_i)=\omega_A(\beta_N(\tilde\lambda_i))$ and
$\omega_e(c_{i})=\omega_{A,1}(\beta_N(\tilde\lambda_i))$ for all $i$.\\
If $\lambda_i$ and $\lambda_j$ were any two parts with odd color,
then $\lambda_i - \lambda_j\ge M+\sum_{\ell = i}^{j-1}
\omega_e(\lambda_{\ell}) \pmod{M}$ translates to if
$\beta_N(\tilde\lambda_i)$ and $\beta_N(\tilde\lambda_j)$ use
$a_1$ in their defining sums, then
$\tilde\lambda_i-\tilde\lambda_j\ge MN+
N(-\omega_A(\beta_N(\tilde\lambda_i))+\omega_A(\beta_N(\tilde\lambda_{j})))
+\beta_N(\tilde\lambda_i)-\beta_N(\tilde\lambda_{j})+
N\sum_{\ell=i}^{j-1}\omega_{A,1}(\tilde\lambda_{\ell})$. \qed\\

\section{Concluding Remarks} \setcounter{equation}{0}
There are undoubtedly many more applications of the iterative
methods described in this paper.  To help motivate the iterative
process, we have described how to reproduce and generalize some
famous results by executing a rearrangement of colors at each
step. However, the Alladi-Gordon bijection may also be ``naively"
iterated without performing this rearrangement.  It would be
worthwhile to investigate the theorems on colored partitions that
arise in this way.  Also, as observed in \cite{Al-Go2}, one can
write down companions to Schur's theorem by reordering the parts
immediately before Step 4 of the proof of Theorem \ref{alladi}
presented in Section 2. Such reorderings could probably be applied
in the case of Theorem \ref{lovejoy1} as well. Finally, it will be
noted that the products in \eqref{product1} and \eqref{product2}
are rather asymmetric.  A more general problem of the nature
considered here would be to give an interpretation of the infinite
product
\begin{equation}
\prod_{k=1}^{\infty} (1+y_1q^{M_1k - R_1})(1+y_2q^{M_2k -
R_2})\cdots(1+y_nq^{M_n k - R_n})
\end{equation}
in terms of partitions whose parts occur in $2^{n}-1$ colors and
satisfy some tractable difference conditions.

\section*{Acknowledgments} The authors are pleased to thank the
Center for Combinatorics at Nankai University, where this work was
initiated, and the Center of Excellence in Complex Systems and the
Department of Mathematics and Statistics at the University of
Melbourne, where this work was completed.

\end{document}